\documentclass[12pt]{amsart}
\usepackage[all]{xy}
\usepackage[dvips]{graphicx}

\newtheorem{teo}{Theorem}[section]
\newtheorem{defi}[teo]{Definition}
\newtheorem{lema}[teo]{Lemma}
\newtheorem{obs}[teo]{Remark}
\newtheorem{cor}[teo]{Corollary}

\newcommand{\bdefi}{\begin{defi}}
\newcommand{\edefi}{\end{defi}}
\newcommand{\bteo}{\begin{teo}}
\newcommand{\eteo}{\end{teo}}
\newcommand{\blem}{\begin{lema}}
\newcommand{\elem}{\end{lema}}
\newcommand{\bobs}{\begin{obs}}
\newcommand{\eobs}{\end{obs}}
\newcommand{\bcor}{\begin{cor}}
\newcommand{\ecor}{\end{cor}}

\newcommand{\Or}{\mathcal{O}}
\newcommand{\fimdem}{\phantom.\hfill \square}
\newcommand{\Z}{\mathbb{Z}}
\newcommand{\dem}{{\bf Proof: }}
\newcommand{\ov}{\overline}

\newcommand{\wht}{\widehat}

\begin{document}

\title[Jiang-type theorems for coincidences]{Jiang-type theorems for coincidences of maps into homogeneous spaces}
\author{Daniel Vendr\'uscolo}
\address{Departamento de Matem\'atica - UFSCar, S\~ao Carlos, SP - Brazil.}
\email{daniel@dm.ufscar.br}

\author{Peter Wong}
\address{Department of Mathematics, Bates College. Lewiston, ME 04240 - USA}
\email{pwong@bates.edu}
\date{\today}
\thanks{This work was initiated during the first author's visit to Bates College, April 13 - May 16, 2006 and was completed during the second author's visit to S\~ao Carlos, November 18 - 23, 2006. The authors would
like to acknowledge support from FAPESP and the National Science Foundation, respectively.}

\begin{abstract}
Let $f,g: X\to G/K$ be maps from a closed connected orientable manifold $X$ to an orientable coset space $M=G/K$ where $G$ is a compact connected Lie group, $K$ a closed subgroup and $\dim X=\dim M$. In this paper, we show that if $L(f,g)=0$ then $N(f,g)=0$; if $L(f,g)\ne 0$ then $N(f,g)=R(f,g)$ where $L(f,g), N(f,g)$, and $R(f,g)$ denote the Lefschetz, Nielsen, and Reidemeister coincidence numbers of $f$ and $g$, respectively. When $\dim X> \dim M$, we give conditions under which $N(f,g)=0$ implies $f$ and $g$ are deformable to be coincidence free.
\end{abstract}


\maketitle

\section{Introduction}

In classical Nielsen fixed point theory, the Nielsen number is either zero or is equal to the Reidemeister number if the so-called Jiang condition is satisfied. For coincidences, we say that a closed connected orientable manifold $M$ is of {\it Jiang-type for coincidences} if for any closed connected orientable manifold $X$ with $\dim X=\dim M$ and for any maps $f,g:X\to M$ we have $L(f,g)=0 \Rightarrow N(f,g)=0$ or $L(f,g)\ne 0 \Rightarrow N(f,g)=R(f,g)$ where $L(f,g), N(f,g)$, and $R(f,g)$ denote the Lefschetz, Nielsen, and Reidemeister coincidence numbers respectively. Evidently, classical Jiang spaces are of Jiang-type for coincidences. In addition, certain $\mathcal C$-nilpotent spaces for the class $\mathcal C$ of all finite groups \cite{dape1} which include coset spaces $G/K$ of compact connected Lie groups $G$ by finite subgroups $K$ and compact nilmanifolds \cite{dape2} are known to be of Jiang-type for coincidences. Recently, it was shown \cite{pe3} that orientable!
  coset spaces $G/K$ of compact connected Lie groups are of Jiang-type (for fixed points). The main purpose of this paper is to show that orientable coset spaces are indeed of Jiang-type for coincidences. We also investigate certain situations when the dimension of the domain space is greater than that of the target space. We assume that the reader is familiar with the geometric Nielsen coincidence theory and root theory as developed in \cite{rob2}.

\section{Coincidences of maps into homogeneous spaces} \label{gs}

Let $G$ be a compact, Lie group and $K$
a closed subgroup of $G$. We will assume the homogeneous space
 of left cosets $G/K$ is orientable. We denote by $p:G\to G/K$ the
natural projection. In general the projection $p$ is a
locally trivial fibration (see \cite[page 675]{wh} and
\cite[\S 7]{st}) with fibre $K$ and in case that $K$ is finite it is
a covering map. The group $G$ acts on $G/K$ by $g\cdot g'K=gg'K$ and $K$ acts on $G$ by $k\cdot g=gk^{-1}$. We should point out that the map $p$ is not a $K$-map.

Let $X$ be an orientable
manifold and $f,h:X\to G/K$ maps. We define $q:\wht X\to X$
to be the induced bundle over $X$ by $f$
i.e., $\wht X=\{(x,g)\in X\times G\ |\ f(x)=p(g)\}$ and
$q(x,g)=x$. By \cite[\S 10]{st}, $q$ is a locally trivial
fibration and $\wht X$ is a $K$-space where the $K$ action on $\wht X$ is the diagonal action on $X\times G$ and $K$ acts trivially on $X$. We have the following
commutative diagram.

$$\xymatrix{
\wht{X} \ar[r]^-{\wht{f}} \ar[d]_-{q} & G \ar[d]^-{p} \\
X \ar[r]^-{f} & {G/K}}
$$

Now we define a map $\eta:\wht X\to G/K$ by:
$$\eta(\wht{x})=[\wht{f}(\wht x)]^{-1}\cdot h(q(\wht x)).$$

\blem \label{coinroot} A point $x \in X$ is a coincidence of the pair $(f,h)$ iff for each
$\wht x \in q^{-1}(x)$, $\eta(\wht x)=eK$.
\elem

\dem ($\Rightarrow$) Suppose $x\in Coin(f,h)$, i.e. $f(x)=h(x)$ , then for each
$\wht x\in q^{-1}(x)$ we have:
\begin{align}
\eta(\wht x)&= [\wht{f}(\wht x)]^{-1}\cdot h(q(\wht x)) \nonumber \\
&= [\wht{f}(\wht x)]^{-1}\cdot h(x) \nonumber \\
&= [\wht{f}(\wht x)]^{-1}\cdot f(x) \nonumber \\
&= [\wht{f}(\wht x)]^{-1}\cdot \wht{f}(\wht x)K \nonumber \\
&= eK. \notag
\end{align}

($\Leftarrow$)
\begin{align}
\eta(\wht x)= eK &\Rightarrow (\wht{f}(\wht x))^{-1}\cdot h(q(\wht x))=eK \nonumber \\
& \Rightarrow [\wht{f}(\wht x)]^{-1}\cdot h(x)=eK \nonumber \\
&\Rightarrow \wht{f}(\wht x)K= h(x) \notag \\
&\Rightarrow f(x)=h(x). \notag
\end{align}
$\fimdem$

\blem \label{downclass} If $\wht{c_1}, \wht{c_2}\in \wht{X}$ are Nielsen equivalent as roots of $\eta$ then
their projections $c_1, c_2 \in Coin(f,h)$ are  Nielsen equivalent as coincidences of $f$ and $h$.
\elem

\dem Let $\wht{c_1}$ and $\wht{c_2}$ be two roots of $\eta$ in the same
Nielsen root class. Then there exists a path $\wht{\gamma}$ from $\wht{c_1}$
and $\wht{c_2}$ such that $\eta\circ\wht{\gamma}$ is homotopic to the trivial
loop at $eK$ relative to the endpoints so that there exists a homotopy
$\wht{H}:I\times I\to M$ such that:
$$\wht{H}(t,0)=\eta(\wht{\gamma}(t))=[\wht{f}(\wht{\gamma}(t))]^{-1}\cdot h(q(\wht{\gamma}(t))),
\ \ \ \wht{H}(t,1)=eK,$$
and
$$\wht{H}(0,s)=\wht{H}(1,s)=eK.$$

The projection $\gamma=q(\wht{\gamma})$ is a path from $c_1$ to $c_2$ in $X$,
 and $H(t,s)=[\wht{f}(\wht{\gamma}(t))]\cdot \wht{H}(t,s)$ is a homotopy between $h\circ\gamma$
and $f\circ \gamma$ relative to endpoints. This implies that
 $c_1$ and $c_2$ are in the same Nielsen coincidence class.
$\fimdem$

\blem \label{upclass} If $c_1, c_2\in Coin(f,h)$ are Nielsen equivalent then for each
$\wht{c_1}\in q^{-1}(c_1)$ there exists $\wht{c_2}\in q^{-1}(c_2)$ that is
Nielsen equivalent to $\wht{c_1}$ as roots of $\eta$.
\elem

\dem If $\gamma$ is a path from $c_1$ to $c_2$ such that
$f\circ \gamma\sim h\circ\gamma$ relative to endpoints, one can lift $\gamma$
to a path $\wht{\gamma}$ starting at $\wht{c_1}$ so that $\wht{\gamma}(1)\in q^{-1}(c_2)$
and
\begin{align}
\eta(\wht{\gamma}(t))&=[\wht{f}(\wht{\gamma}(t))]^{-1}\cdot h(q(\wht{\gamma}(t))) \notag \\
&= [\wht{f}(\wht{\gamma}(t))]^{-1}\cdot h(\gamma(t)) \notag \\
&\sim [\wht{f}(\wht{\gamma}(t))]^{-1}\cdot f(\gamma(t)) \notag \\
&= eK. \notag
\end{align}
Note that the homotopy used here is relative to the endpoints.
$\fimdem$

\bobs These two lemmas show that there is an injection from the (non-empty) coincidence classes
of $(f,h)$ to the (non-empty) roots classes of $\eta$ and a surjection from the
(non-empty) roots classes of $\eta$ to the (non-empty) coincidence classes
of $(f,h)$.
\eobs

Recall that the space $\wht{X}$ is a $K$-space and the action of $K$ in $\wht{X}$ is trivial on the first coordinate. From \cite{rob2}, an (geometrically) inessential coincidence class is a class that disappears under homotopy.

\bteo \label{ines}\begin{enumerate}
\item If a coincidence class $C\subset Coin(f,h)$ is (geometrically) inessential then all root classes of $\eta$ associated to $C$ by lemmas~\ref{upclass} and  \ref{downclass} are all (geometrically) inessential by a $K$-homotopy.
\item If the map $\eta$ can be made root free by a $K$-homotopy then the pair $(f,h)$ can be made coincidence free.
\end{enumerate}
\eteo

\dem (1) Suppose that $F, H:X\times I\to G/K$ are homotopies such that:
$$F(x,0)=f(x),\ \ F(x,1)=f'(x),$$
$$H(x,0)=h(x),\ \ H(x,1)=h'(x)$$
and $C$ is not related to any coincidences of $h'$ and $f'$. Since $G/K$ is a manifold, it follows that one can assume the homotopy $F$ constant so that $F(x,t)=f(x)$ for all $t$.
Let $\wht{F}:\wht{X}\times I\to G$ be the constant homotopy with
$\wht{F}(\wht{x},t)=\wht{f}(\wht x)$.

Now we define $\wht{\mathcal F}:\wht{X}\times I\to G/K$ by:
$$\wht{\mathcal F}(\wht{x},t)=[\wht{F}(\wht{x},t)]^{-1}\cdot H(q(\wht{x}),t)).$$
The map $\wht{\mathcal F}$ is a $K$-homotopy and we have:
$$\wht{\mathcal F}(\wht{x},0)=[\wht{F}(\wht{x},0)]^{-1}\cdot H(q(\wht{x}),0))=[\wht{f}(\wht{x})]^{-1}\cdot h(q(\wht{x}))=\eta(\wht{x}),$$
and
$$\wht{\mathcal F}(\wht{x},1)=[\wht{F}(\wht{x},1)]^{-1}\cdot H(q(\wht{x}),1))=[\wht{f}(\wht{x})]^{-1}\cdot h'(q(\wht{x})).$$

Now all root classes of $\eta$ associated to the coincidence class $C$ in the
sense of lemmas~\ref{upclass} and \ref{downclass} will be inessential by
the $K$-homotopy $\wht{\mathcal F}$, that is, the the root classes of $\eta$ corresponding to $C$ are not related to any root class of $\wht {\mathcal F_1}$.

(2) Let $\wht{\mathcal F}:\wht{X}\times I\to G/K$ be a $K$-homotopy such that $\wht{\mathcal F}(0,t)=\eta$. We define $F:X\to G/K$ by:
$$F(x,t)=[\wht{f}(\wht{x})]\cdot \wht{\mathcal F}(\wht{x},t)\ \ \ \mbox{where}\ \ q(\wht{x})=x.$$

If $q(\wht{x})=q(\wht{x'})$ then there exists $k\in K$ such that $\wht{x}=(x,g)$ and $\wht{x'}=(x,gk)$ for some $g\in G$ and $k\in K$. The homotopy $\wht{\mathcal F}$ is $K$-equivariant and the diagram
$$\xymatrix{
\wht{X} \ar[r]^-{\wht{f}} \ar[d]_-{q} & G \ar[d]^-{p} \\
X \ar[r]^-{f} & {G/K}}
$$
commutes. It follows that $
[\wht f(\wht x')]\cdot{\wht {\mathcal F}}(\wht x',t)=[\wht f(\wht x)]\cdot{\wht {\mathcal F}}(\wht x,t)
$ so that the homotopy $F$ is well defined.

We note that:
$$F(x,0)=[\wht{f}(\wht{x})]\cdot \wht{\mathcal F}(\wht{x},0)=[\wht{f}(\wht{x})]\cdot\eta(\wht{x})=[\wht{f}(\wht{x})]\cdot [\wht{f}(\wht{x})]^{-1}\cdot h(q(\wht{x}))=h(x),$$
and if $\wht{\mathcal F}(\wht{x},1)^{-1}(eK)=\emptyset$ then we define $h':X\to G/K$ by:
$$h'(x)={F}(x,1)=[\wht{f}(x)]\cdot \wht{\mathcal F}(\wht{x},1)\neq f(x)\ \ \forall x\in X$$
so that we have $h\sim h'$ and $Coin(f,h')=\emptyset$. $\fimdem$


Thus, we have just shown that the Nielsen coincidence problem for the pair $(f,h)$ is equivalent to the $K$-equivariant Nielsen root problem of $\eta$. More precisely, we have the following result.

\bcor \label{khomo} Let $X,M$ be manifolds. Suppose $M=G/K$ is the homogeneous space of left cosets of a Lie group $G$ by a closed subgroup $K$. For any maps $f,h:X\to M$ we denote $q:\wht{X}\to X$ the induced fibration over $X$ of the fibration $p:G\to G/K$ by the map $f$ and we define the map $\eta:\wht{X}\to G/K$ by
$\eta(\wht{x})=[f(\wht x)]^{-1}\cdot h(q(\wht x))$. Then there exists a pair $(f',h')\sim (f,h)$ such that $Coin(f',h')=\emptyset$ iff the map $\eta $ can be made root free by a $K$-homotopy in $\wht{X}$.
\ecor

\section{Main results - codimension zero}

In this section we assume that $G/K$ and $X$ are connected, compact, orientable
$n$-manifolds. For $n\ge 3$ these spaces are Wecken spaces for coincidences and
if $n=1$ then $G/K=S^1$, for $n=2$, $G/K$ is the Torus or $S^2$ (all such spaces are Wecken spaces for coincidences).
In all cases we suppose that the coincidence classes of $(f,h)$ are
finite, then $Coin(f,h)=\{c_1,\dots,c_r\}$ or empty.
We will follow the approach of \cite{pe0} and \cite{pe3}.

For every coincidence point $c_i$, each point of the fibre $q^{-1}(c_i)$
is a root and this set has the same number $m$ of connected components as $K$.
Moreover, $\eta$ has no other roots than these fibres.

It is easy to see that each connected component of $q^{-1}(c_i)$ are in the same root class of $\eta$ and each root class of $\eta$ is the union of some connected components of $q^{-1}(c_{i_1},c_{i_2},\cdots,c_{i_l})$ ($\{c_{i_1},c_{i_2},\cdots,c_{i_l}\}$ a coincidence class of $(f,h)$)

Denote by $\Or_i^j$ ($1\le i\le r$, $1\le j\le m$ and $q(\Or_i^j)=c_i$)
each connected component of the root set, and $\Or_i=q^{-1}(c_i)$.
The fibration $q:\wht{X}\to X$ is locally trivial and we can choose
neighborhoods $V_i\subset X$ of $c_i$ such that $\ov{V_i}\cap Coin(f,h)=\{c_i\}$,
$q^{-1}(V_i)=\bigcup\limits_{j=1}^{m}V_i^j\times K$
with $\Or_i^j\subset V_i^j$,  then we can see
that $(q|_{V_i^j})_{*}:H_n(V_i^j,V_i^j\setminus \Or_i^j)\to H_n(V_i,V_i\setminus c_i)\cong \Z$
is an isomorphism.

 Let $o_i\in H_n(V_i,V_i\setminus c_i)$ be the fundamental class around the coincidence point $c_i$ and
$\wht{o_i^j}=(q|_{V_i^j})_{*}^{-1}(o_i)$ (the ``fundamental class"
 around $\Or_i^j$) and $\wht{o_i}= \bigoplus\limits_{j=1}^m \wht{o_i^j}\in H_n(q^{-1}(V_i),q^{-1}(V_i)\setminus \Or_i) $ (the ``fundamental class" around $\Or_i$).

Define the
{\it local root index} of $\Or_i$ with respect to $\eta$ to be the integer $\omega_K(\eta;\Or_i)$ that is the image of $\wht{o_i}$ by the map:
$$H_n(q^{-1}(V_i),q^{-1}(V_i)\setminus \Or_i)\stackrel{\eta_*}{\longrightarrow} H_n(M,M\setminus eK)\cong H_n(M)\cong \Z$$

Similarly for each $j$, we define $\omega_K(\eta;\Or_i^j)$ by the image of $\wht{o_i^j}$ by the restriction of $\eta_*$ to $H_n(V_i^j, V_i^j\setminus \Or_i^j)$.

Since $\eta$ is $K$-equivariant, we have
$\omega(\eta; \Or_i^{j_1})=\omega_K(\eta; \Or_i^{j_2})$ for $1\le j_1,j_2\le m$.
Then we conclude $\omega_K(\eta;\Or_i)=m\cdot\omega(\eta; \Or_i^{1})$.

Since $q^{-1}(c_i)$ has the same number, $m$, of connected components of $K$, we denote by $\Or_i^1$ the one associated with the connected component of $K$ that contains $e\in G$.

\blem \label{equal} The local coincidence index of $c_i$ with respect to $(f,h)$ coincides with the local root index of $\Or_i^1$ with respect to $\eta$, i.e., $ind(f,h;c_i)=\omega_K(\eta;\Or_i^1)$.
\elem

\dem Let $\sigma:(\Delta^n,\partial\Delta^n)\to (X, X\setminus c_i)$ be a singular $n$-simplex such that $[\sigma]=o_i\in H_n(X, X\setminus c_i)$ and let $\wht{\sigma}: (\Delta^n,\partial\Delta^n)\to (\wht{X}, \wht{X}\setminus \Or_i^1)$ be a lift of $\sigma$ i.e. $q\wht{\sigma}=\sigma$. Then $[\wht{\sigma}]=o_i^1\in  H_n(\wht{X}, \wht{X}\setminus \Or_i^1)$

We know $ind(f,h;c_i)=[(f\sigma,h\sigma)]\in H_n(M\times M, M\times M\setminus \Delta M)$ and $\omega_K(\eta;\Or_i^1)=[(\wht{f}\wht{\sigma})^{-1}\cdot h\sigma] \in H_n(M,M\setminus eK)$.

Now if $\tau_t$ is a homotopy from $(\wht{f}\wht{\sigma})^{-1}$ to the constant map at $eK$, we can see that the homotopy $[(\tau_f\cdot f\sigma, \tau_t\cdot h\sigma)]$ implies that $[(f\sigma,h\sigma)]=[eK,(\wht{f}\wht{\sigma})^{-1}\cdot h\sigma)]$.

Then we just note that the map $\xi:(M,M\setminus eK)\to (M\times M, M\times M\setminus \Delta M)$ defined by $\xi(x)=(eK,x)$ induces an isomorphism on homology, therefore $ind(f,h;c_i)=\omega_K(\eta;\Or_i^1)$. $\fimdem$

For any $i,j$, $\Or_i^j$ is homeomorphic to $\Or_i^1$ by an orientation preserving homeomorphism which preserves the root index. Thus, we have

\bcor $ind(f,h;c_i)=\omega_K(\eta;\Or_i^j),\ \ \forall j\  1\le j \le m.$
\ecor

\blem \label{samesig} If $\omega_K(\eta;\Or_i^j)=0$ for some $j$, then $\omega_K(\eta;\Or_i^l)=0$ for all $1\le i\le r$ and $1\le l\le m$. If $\omega_K(\eta;\Or_i^j)\neq 0$, then
$$\omega_K(\eta;\Or_i^j)\cdot \omega_K(\eta;\Or_s^l)>0\mbox{  for all }
1\le s\le r\mbox{ and } 1\le l\le m.$$
\elem

\dem It is straightforward to verify that the proof presented in \cite[Lemma 2.1]{pe3}
remains valid under our hypothesis if we replace
 the map $\varphi$ there  by our map $\eta$. $\fimdem$

\bcor \label{nonemp} All non-empty coincidence classes of the pair $(f,h)$ have either index zero or index of the same sign.
\ecor

\dem It follows from Lemmas~\ref{equal} and \ref{samesig}. $\fimdem$

Our main result in this section is that $L(f,h)\neq 0\Rightarrow N(f,h)=R(f,h)$. We know that for each Reidemeister coincidence class we can associate one (possible empty) Nielsen coincidence class. This correspondence depends on many choices (base points, lifts, etc.). We established that all non-empty classes have index of same sign and the next lemma shows that the inessential classes can be considered as non-empty.

\blem \label{reidclass} Given $\beta\in\pi_1(M)$ and a pair of maps $f,h:X\to M$ such that $X$ and $M$ are connected manifolds, there exists a map $h'\sim h$ and a point $x_{1}\in Coin(f,h')$ so that the coincidence class to which $x_{1}$ belongs is associated to the Reidemeister class of $\beta$.
\elem

\dem Suppose that $Coin(f,h)\neq \emptyset$ and $Coin(f,h)\neq M$. Let $x_0\in Coin(f,h)$ and choose $x_{1}\in (M\setminus Coin(f,h))$ and a path $\gamma:[0,1]\to X$ such that $\gamma(0)=x_0$ and $\gamma(1)=x_1$.

Let $B_{h(x_1)}$ be an open neighborhood of $h(x_1)$ such that $f(x_1)\notin \ov{B_{h(x_1)}}$. The connected component $O$ of $h^{-1}(B_{h(x_1)})$ that contains $x_1$ is an open set and there exists an open neighborhood $B_{x_1}$ of $x_1$, such that:
\begin{enumerate}
\item $\ov{B_{x_1}}\subset O$;
\item $\ov{B_{x_1}}$ is homeomorphic under $\psi$ to the unit ball $B(0,1)$ (of same dimension of $X$) such that $\psi(x_1)=0$.
\item $\ov{B_{x_1}}\cup Coin(f,h)=\emptyset $ (since $Coin(f,h)$ is a closed set).
\end{enumerate}

For each point $x\in \ov{B_{x_1}}$ define $x_{\partial}\in \partial B(0,1)$ as the image of $\psi(x)$ by the radial projection, then there exists a unique $t_x\in [0,1]$ such that $\psi(x)=t_x\cdot x_{\partial}$.

Now we define $h':X\to M$ by
$$h'(x)=\left\{
\begin{array}{lcl}
h(x)                                           & if & x\in X\setminus B_{x_1}\\
h(\psi^{-1}((1-4t_x)\cdot x_{\partial}))       & if & 0\le t_x\le 1/4\\
h(\gamma(2-4t_x))                              & if & 1/4 \le t_x\le 1/2\\
\beta(4t_x-2)                                  & if & 1/2 \le t_x\le 3/4\\
f(\gamma(4t_x-3))                              & if & 3/4 \le t_x\le 1.\\
\end{array}
\right.
$$

It is easy to see that $h'$ is the desired map.

If $Coin(f,h)\neq \emptyset$ then a similar homotopy can produce a coincidence point and if $Coin(f,h)=M$ then a local small deformation of one of the maps is enough. $\fimdem$

\bteo \label{main} Let $G$ be a compact connected Lie group, $K$ a closed subgroup and $M=G/K$ the homogeneous space of right cosets. Assume that $M$ and $X$ are orientable n-manifolds and that $X$ is compact. For any maps $f,h:X\to M$ we have:
\begin{enumerate}
\item $L(f,h)=0\Rightarrow N(f,h)=0$
\item $L(f,h)\neq 0\Rightarrow N(f,h)=R(f,h).$
\end{enumerate}
\eteo

\dem Lemma~\ref{reidclass} and Corollary~\ref{nonemp} show that if one Reidemeister class is associated with an empty Nielsen class then all Nielsen classes have index zero (it is simple to see that the coincidence produced in Lemma~\ref{reidclass} has index zero).
Now we complete the proof by recalling that $L(f,h)=\sum ind(f,h;C_i)$. $\fimdem$

\bobs In \cite{dape1}, it was shown that a closed orientable $\mathcal C$-nilpotent manifold whose fundamental group has a center of finite index is a Jiang-type space for coincidences. While orientable coset spaces $G/K$ may well be such a space but we are unable to verify at this point. Moreover, our approach of studying equivariant roots of the $K$-map $\eta$ has the advantage of allowing us to study positive codiemnsional coincidence problem as we illustrate in the next section whereas the approach of \cite{dape1} is restricted to codimension zero situation.
\eobs

\section{Positive codimension}

In this section the Lie group $G$ is not necessarily compact and the homogeneous
space $G/K$ and the manifold $X$ can be non-orientable. In this situation, the classical Lefschetz coincidence number is not well-defined while the Nielsen coincidence number can be defined using the geometric concept of essentiality of \cite{rob2}.

\bteo \label{liftf} Let $X, M$ be manifolds.
Suppose $M=G/K$ is the homogeneous space of left cosets of a connected
Lie group by a closed subgroup $K$ and X is compact.
Let $f,h:X\to M$ be maps. Suppose there exists a map $\wht{f}$ such that the following diagram is commutative.
$$\xymatrix{
& G \ar[d]^-{p} \\
X \ar[ur]^-{\wht{f}} \ar[r]^-{f} & {G/K}}
$$
If $N(f,h)=0$ (the ``geometric" Nielsen number) then there
exist a map $h'\sim h$ such that $Coin(f,h')=\emptyset$
\eteo

\dem First, we define $\eta:X\to G/K$ by $\eta(x)=[\wht{f}(x)]^{-1}\cdot h(x)$ and
we note that Lemmas~\ref{coinroot}, \ref{upclass}, \ref{downclass} and Theorem~\ref{ines} hold true for this map $\eta$ (and there is no action of $K$).

By Lemma~\ref{ines} and \cite{rob1}, $N(f,h)=0 \Rightarrow N(\eta;eK)=0$ (the Nielsen root number).

Now we use the results of \cite{dape3} to obtain a homotopy $\wht{H}:X\times I\to M$
such that :
$$\wht{H}(x,0)=\eta(x),\ \ \mbox{and}\ \ \wht{H}(x,1)\neq eK\ \ \forall x\in X,$$
it means that the map $\wht{H}(x,1):X\to M$ is root free with respect to the point
$eK\in M$. Now we use Corollary~\ref{khomo}.$\fimdem$

\bobs If $f$ is the projection from $G$ to $G/K$, the first part of the
hypothesis of Theorem~\ref{liftf} is satisfied (one can lift the map $f$ putting $\wht{f}=id$).
If we study self-coincidences of the projection, the map $\eta$ will be the
constant map. If $N(\eta;eK)=0$ then we can conclude that there exist
$p'\sim p$ such that $Coin(p,p')=\emptyset$ (this case is covered by
\cite[Proposition 2.3]{dado}). This suggests a different approach to the self-coincidence problem studied in \cite{dado} in which primary and higher obstructions are the principal tools. This will be the objective of our future work.
\eobs

\bobs The existence of $\wht f$ in the hypothesis of Theorem \ref{liftf} can be determined using classical obstruction theory to lifting as in \cite{wh}.
\eobs

Next, we apply our techniques to coincidences into suitable manifolds.
For a manifold $M$, $e\in M$, and $G(M)$ the group of all homeomorphisms of $M$, we say (see \cite{fanu}) that $M$ is {\it suitable} if there exists a continuous map $\theta:M\to G(M)$ such that:
$$\theta(x)(x)=e\ \ \mbox{ and }\ \  \theta(e)=id_M.$$

Brown \cite{br} showed that on a suitable manifold $M$ there exists a multiplication such that:
\begin{enumerate}
\item $x\cdot e=x,\ \ \forall x\in M;$
\item given $a,b\in M, \ \ \exists x\in M$ such that $a\cdot x=b$;
\item $x\cdot y=x\cdot z\Rightarrow y=z,\ \ \forall x, y, z\in M$.
\end{enumerate}
and that for compact manifolds the existence of this multiplication is equivalent to the definition. Moreover Brown notes that for each $x\in M$ there exists a unique inverse $x^{-1}$ such that $x\cdot x^{-1}=e$.

Let $X$ be a compact manifold, $M$ be a suitable manifold and $f,h:X\to M$ be maps. We define $\eta:X\to M$ by $\eta(x)=(f(x))^{-1}\cdot h(x)$ and it is easy to see the arguments used in section~\ref{gs} can be used under these hypotheses showing that there exists a bijection between the coincidence classes of $(f,h)$ and the root classes of $\eta$. Furthermore, the proof of Theorem~\ref{liftf} remains valid even for a non-compact suitable manifold $M$.

\bteo \label{suit} Let $X, M$ be manifolds such that $X$ is compact and $M$ is suitable. If $f,h:X\to M$ are maps such that $N(f,h)=0$ (the geometric Nielsen number defined in \cite{rob1}) then there exist a map $h'\sim h$ such that $Coin(f,h')=\emptyset$.
\eteo

\bobs We should point out that a Lefschetz type coincidence index homomorphism was developed in \cite{saveliev} for the setting of the positive codimension coincidence problem. In particular, a formula for the Lefschetz coincidence index homomorphism was obtained in \cite[Theorem 7.1]{saveliev} when the target is a suitable manifold. However, the vanishing of this index homomorphism does not guarantee the vanishing of $N(f,h)$.
\eobs

\noindent
{\bf Acknowledgment.} We thank the anonymous referee for his/her remarks that help improve the exposition of the paper, in particular, the proof of Lemma \ref{equal}.

\end{document}